
\input amstex
\documentstyle{amsppt}
\NoRunningHeads
\TagsAsMath
\TagsOnRight
\NoBlackBoxes
\font\special=cmr6
\topmatter
\title
On the deformation of abelian integrals.
\endtitle
\author
Feodor A. Smirnov\\
{\special Laboratoire de Physique Theorique
et Hautes Energies,}\\ {\special
4 place Jussieu,
Paris 75005, FRANCE }\\
{\special and}\\
{\special Steklov Mathematical Institute,}
\\{\special Fontanka 27, St. Petersburg 191011,  RUSSIA}
\endauthor
\abstract
We consider the deformation of abelian integrals which
arose from the study of SG form factors. Besides the
known properties they are shown to satisfy Riemann bilinear
identity. The deformation of intersection number of cycles
on hyperelliptic curve is introduced.
\endabstract
\endtopmatter

The subject of this paper is related to the problem
of description of all local fields in SG theory
the solution to which will be published elsewhere.
However, the result presented in this paper
which is the generalization of Riemann bilinear identity
to deformed abelian integrals
constitutes rather
isolated and, presumably, the most mathematically
interesting
part of the problem, for that reason the author
decided to publish it separately.

The integrals which
provide a deformation of periods of hyperelliptic differentials
were introduced several years ago for the description of
SG form factors [1], let us explain how they are constructed.
Consider the function
\def\a{\alpha}
\def\b{\beta}
$$
\varphi (\a)=C
\ \text{exp}\left(-2\int\limits _0^{\infty}
{{\text{sin}^2{k\a\over 2}\ \text{sh}{\pi+\xi\over 2}k}\over
{k\ \text{sh}{\xi k\over 2}\ \text{sh} \pi k}} dk\right) $$
where $C$ is certain constant [1] which is needed for $\varphi$
to satisfy the relation
$$
\varphi(\a +\pi i)\varphi (\a)=\frac 1 {4\ \text{sh}{\pi\over\xi}
(\a+{\pi i\over 2})
\ \text{ch}(\a)} \tag 1 $$
without any additional multiplier in RHS.
Other important properties of $\varphi$ are $$
\varphi (\a+i\xi)= \varphi (\b)\ \frac {\text{ch}{1\over 2}(\a-{\pi i\over 2})}
{\text{ch}{1\over 2}(\a+{\pi i\over 2}+i\xi)}
,\qquad\varphi (\a+2\pi i)= -\varphi (\a)
\ \frac {\text{sh}{\pi\over \xi}(\a+{\pi i\over 2})}
{\text{sh}{\pi\over \xi}(\a+{3\pi i\over 2})} \tag 2  $$
The function $\varphi (\a)$ does not have other singularities
for $0\le \text{Im}\a\le 2\pi $ but the simple poles at the
points $\a={\pi i \over 2}+i\xi k,\qquad k\ge 0$ (how many
of them are in the strip depends on the value of $\xi$).
Asymptotically when $\a\to\pm\infty$
the function $\varphi (\a)$ behaves as
$$ \varphi (\a)\sim \text{exp}\bigl( -{1\over 2}({\pi\over\xi}+1)|\a |\bigr)$$

The deformation of hyperelliptic abelian integral is defined
as follows [1,2]. Consider $2n$ points $\b _1,\cdots , \b_{2n}$
(analogues of branching points).
Then for two given polynomials $Q(a)$ and $L(A)$ (which can also
depend respectively
on $b_j=e^{{2\pi\over\xi}\b _j}$
and on $B_j=e^{\b_j}$ as on parameters)
we consider the paring $\langle Q(a),L(A)\rangle$ defined by an integral:
$$
\langle Q(a),L(A)\rangle\equiv
\int\limits _{-\infty}^{\infty}\prod\limits _{j=1}^{2n}
\tilde{\varphi} (\a,\b _j)
\ Q(e^{{2\pi\over\xi}\a})L(e^{\a})
\ e^{{2\pi\over\xi}\a}
\ d\a \tag 3
$$
where for further convenience we introduced the notation:
$$ \tilde{\varphi}(\a ,\b)=
\varphi(\a -\b)\ \text{exp}(-{1\over 2}({\pi\over\xi}+1) (\a+\b))$$
The integral is convergent for
$1\le\text{deg}L(A)\le 2n-1$ and $\text{deg}Q(a)\le n-1$. However, as it
is explained in [1,2] the last condition can be relaxed, with
a proper regularization the integral can be defined for
arbitrary polynomial $Q(a)$
in such a way that one can deal with the deformed integral as
with usual one: deform contours for example. For the sake
of completeness we shall describe the regularization in the Appendix.
Moreover, only $2n-1$ of these
polynomials give really different integrals, we shall comment on this
point later. In what follows we shall use equally often the variables
$$\a ,\b _j ;\qquad a=e^{{2\pi\over\xi}\a}, b_j=e^{{2\pi\over\xi}\b _j};
\qquad A=e^{\a}, B_j = e^{\b _j}$$

Let us explain the relation to hyperelliptic integrals.
Consider the limit when $\xi\to\infty$
but the variables
$b_j$ remain finite, which
means that $\b _j$ are getting simultaneously rescaled. Let us
require that $b_1<b_2<\cdots<b_{2n} $
Then it can be shown
that the following asymptotic formula holds
$$\langle  a^p,A^k\rangle\sim \prod\limits _{j=1}^{k}B_j
\ \int\limits _{\gamma _{k}} {a^p\over \sqrt{P(a)}}
da \tag 4 $$
where the polynomial $P(a)$ is
given by
$P(a)=\prod (a -b _j)$, the cycles $\gamma _j$
are those drawn around the branching points $b_{i}$ and
$b_{i+1}$
on the hyperelliptic surface $c^2=P(a)$.
In the papers [3,4] the formula (4) was proven in a little
different context: the limit $\xi\to\infty$ was taken first,
and the $\b _j$ were rescaled.

The reason the asymptotical formula (4) to
exist is hidden in the properties of the function $\tilde{\varphi}(\a,\b)$
for $\xi\to\infty$:
$$\tilde{\varphi}(\a ,\b )\sim
\bigl(A^{-1}\theta (\a-\b)+B^{-1}\theta (\b-\a)\bigr)
\ {1\over\sqrt{(a-b _j)}},$$
the equation (1) and the second equation from (2) turn into
$$\left({1\over\sqrt{(a-b _j)}}\right)^2 = {1\over(a-b _j)},
\qquad {d\over da}{1\over\sqrt{(a-b_j)}}=
-{1\over 2} {1\over(a-b_j)} {1\over\sqrt{(a-b_j)}} $$

Let us emphasize the most interesting property of this deformation of
hyperelliptic integrals: the integration over different cycles
is replaced by the integration over all the real axis with
additional weights $A^k=e^{k\a}$ inserted into the integral. In what
follows we shall call the monomial
$A^{k}$ a deformation of $\gamma _{k}$.
It is natural to consider usual cycles as generators of an
exterior algebra over the ring of integer numbers, the deformed
cycles constitute a natural exterior algebra over the ring of
quasiconstants: symmetric Laurent polynomials in variables
$B_j$.

Let us list the properties of the hyperelliptic
differentials. We do not expect the Riemann surface itself to
allow a deformation, only the globally defined objects:
cycles and differentials with their periods.
That is why we are interested in the differentials of the
first and of the second kind, i.e. in those without simple poles
because including the differentials with the simple poles is
the same as considering points on the surface.
It is convenient for our goals to consider the following
basis of the differentials of the first and of the second
kind respectively
$$
\align
&\eta _p={a ^{p-1}\over \sqrt{P(a)}}da,\tag 5 \\
&\zeta _p= { 1 \over \sqrt{P(a)}}   \bigl(\sum\limits _{k=n}
^{2n-p} +{1\over 2} \sum\limits _{k=1}^{n-1}\bigr)
\bigr((-1)^{k+p}(k-p)a ^{k-1}\sigma _{2n-p-k}(b_1,\cdots ,b_{2n}
\bigr)da
\endalign $$
where $p=1,\cdots ,n-1$ (recall that the genus of the surface is $n-1$),
$\sigma _{i}(b_1,\cdots ,b_{2n})$ are the elementary
symmetric polynomials.
We are missing the differential with the
polynomial in front of ${ 1 \over \sqrt{P(a)}}$
whose leading power is
$n-1$ in (5) because it is
a differential with the simple pole at the infinity.
The differentials
with the polynomials of higher degree than $2n-2$ can be reduced
to those with lower degree by subtracting the exact forms
$d(a ^k \sqrt{P(a)})$.
We are going to generalize to the deformed case
the following properties of
the differentials and cycles:

{\bf 1. Total derivatives.}
The differentials are defined up to the adding of an exact form.

{\bf 2. The independent cycles.}
There are only $2n-2$ independent cycles, there is one relation of
linear dependence between the cycles $\gamma _i$, namely,
$\sum \gamma _{2i-1}\sim 0$
which means that for any differential without residues on the
surface one has
$$\sum\limits _{i=1}^{n}\int _{\gamma _{2i-1}}\omega =0 \tag 6 $$

{\bf 3. Riemann bilinear identity.} For any two differentials
without residues one has
$$\sum\limits _{i=1}^{n-1} \left(
\int _{\a _i}\omega _1\int _{\b _i}\omega _2
 - \int _{\a _i}\omega _1\int _{\b _i}\omega _2\right)
=\sum\limits _{\text{poles}}\text{res}(\Omega _1 \omega _2)
\equiv\omega _1\circ\omega _2
$$
where $ \Omega _1 $ is a primitive function for  $\omega _1$;
$\a _i,\b _i$ is a canonic basis with the intersection numbers:
$$\a _i\circ \a _j=\b _i\circ \b _j=0,\ \a _i\circ\b _j=\delta _{ij}$$
It can be chosen as $\a _i=\gamma _{2i}$, $\b _i=\sum _{k=1}^{i}\gamma
_{2k-1}$.  The usual way of proving of the relation is to integrate the 1-form
$\Omega _1 \omega _2 $ over the canonic polygon associated to the Riemann
surface. This way is not easy to generalize to the deformed case, so we shall
rely on an alternative possibility which is described later.

The differentials (5) are constructed in order that they
constitute a canonic basis of differentials:
$$ \eta _i\circ\eta _j =  \zeta _i\circ\zeta _j =0,
\ \eta _i\circ\zeta _j=\delta _{ij} $$
This property of the differentials allows to write a dual form
of the Riemann bilinear identity (equivalent to the original one):
$$\sum\limits _{i=1}^{n-1}
\left(\int _{\delta _1}\eta _i
\int _{\delta _2}\zeta _i -
\int _{\delta _2}\eta _i \int _{\delta _1}  \zeta _i \right)
=\delta _1\circ \delta _2 \tag 7 $$
for any two cycles $\delta _1,\delta _2$ on the surface.
This is the formulation which allows easier generalization
for deformed case, so let us sketch the proof of (7).

Consider the identity which follows from a simple
algebra:
$$
 \sum\limits _{i=1}^{n-1} \bigl(\eta _i (a)\zeta _i (a') -
\zeta _i(a)  \eta _i(a')\bigr)=
\bigl({\partial\over \partial a}{\sqrt{P(a)}\over
(a -a')\sqrt{P(a')}}-
{\partial\over \partial a'}{\sqrt{P(a')}\over
(a' -a)\sqrt{P(a)}}\bigr)da da'  $$
Integrate this identity
over $\delta _1\times \delta _2$. In the LHS one gets
the LHS of (7) while in the RHS one has the integral of total
derivatives which is evidently sitting on the contact terms,
the calculation of the latter gives the intersection number of the cycles.

Let us explain how to deform all these properties of the
abelian integrals. The generalization of the first and the
second properties are well known [2,3], but we are going to
discuss them for the sake of completeness.

{\bf 1. The deformation of the total derivative.}
For the given polynomial $Q(a)$ let us construct
the polynomials of the form
$$\widehat{Q}(a)=
\bigl(\prod\limits _{j=1}^{2n} (a\tau-b_j)Q(a\tau ^4)-
\prod\limits _{j=1}^{2n} (a\tau ^{-1}-b_j)Q(a)\bigr)a^{-1}$$
here and later on $\tau=\text{exp}({\pi ^2i\over\xi})$.
We claim that substituting $\widehat{Q}(a)$ into the integral (3)
one gets zero. Formally it goes as follows:
$$ \align
&\langle  \widehat{Q}(a),L(A)\rangle=
\int\limits _{-\infty}^{\infty}\prod\limits _{j=1}^{2n}
\tilde{\varphi}(\a ,\b _j)
\widehat{Q}(e^{{2\pi\over\xi}\a})L(e^{\a})
e^{{2\pi\over\xi}\a}d\a =  \tag 8 \\&=
\left\{\int\limits _{-\infty +2\pi i}^{\infty +2\pi i}
-\int\limits _{-\infty}^{\infty}\right\}
\prod\limits _{j=1}^{2n}
\tilde{\varphi}(\a ,\b _j)(e^{{2\pi\over\xi}\a}\tau ^{-1}
-e^{{2\pi\over\xi}\b _j} )
Q(e^{{2\pi\over\xi}\a})L(e^{\a})
d\a=0
\endalign $$
We used the properties of the function $\varphi (\a)$ to write
this identity. The latter integral is zero because the integrand
does not have singularities in the strip $0<\text{Im}\a <2\pi $.
In fact, both integrals in (8) can be divergent, but the manipulations
with them can be justified with the regularization described at the
Appendix. Notice that in the limit $\xi\to\infty$ the polynomial
$\widehat{Q}(a)$ turns into
$${2\pi\over \xi}\bigl({d\over da}Q(a)-{1\over 2}{d\over da}P(a)\bigr)$$
as it has to be for the total derivative.

{\bf 2. Independent deformed cycles.}
Using the properties of the function $\varphi (\a)$ one can write an
identity
$$ \align
&\int\limits _{-\infty}^{\infty}\prod\limits _{j=1}^{2n}
\tilde{\varphi}(\a ,\b _j)
\bigl\{\prod\limits _{j=1}^{2n}(e^{\a}+ie^{\b _j})-
\prod\limits _{j=1}^{2n}(e^{\a}-ie^{\b _j})\bigr\}
Q(e^{{2\pi\over\xi}\a})e^{{2\pi\over\xi}\a}
d\a =
\\&=
\left\{\int\limits _{-\infty +i\xi}^{\infty +i\xi}
-\int\limits _{-\infty}^{\infty}\right\}
\prod\limits _{j=1}^{2n}
\tilde{\varphi}(\a ,\b _j)(e^{\a}-ie^{\b _j})
Q(e^{{2\pi\over\xi}\a})e^{{2\pi\over\xi}\a} d\a  \tag 9
\endalign $$
Formally, the integral in the RHS is zero (the integrand is regular)
unless we pick up the contributions from the segments
$(-\infty ,-\infty +i\xi)$ and
$(\infty ,\infty +i\xi)$. The latter
happens for the analogues of the third kind differential with
a simple pole at infinity.
Again, one has to be careful with possible divergences and to
make the regularization (see Appendix). For the analogues
of the differentials of the first an of the second kind
one shows that the
integral is zero, hence, decomposing the LHS of (9) one has:
$$
\sum\limits _{k=0}^{n-1} \sigma _{2n-2k-1}(B_1,\cdots ,B_{2n})
\langle Q(a),A^{2k+1}\rangle=0 $$
This identity is analogous to
(6), the only new point being that the deformed cycles are
linear dependent over the ring of symmetric functions
of $B_j$ (quasiconstants).

It is wonderful that the properties 1 and 2 being of very
different nature in underformed case allow quite similar
proofs after the deformation. That means that cycles
and differentials become much more similar after the deformation.

{\bf The deformation of Riemann bilinear identity}.
Let us introduce the basis of deformed first and second kind differentials
similar to (5). The corresponding polynomials are given by
$$\align
& R_p(a)=a^{p-1},\\ &
S_p (a)=\bigl(\sum\limits _{k=n}^{2n-p}
+{1\over 2} \sum\limits _{k=1}^{n-1}\bigr)
\bigr((-1)^{k+p}
(\tau ^{k-p} -\tau ^{p-k})
a ^{k-1}\sigma _{2n-p-k}(b_1,\cdots ,b_{2n})
\bigr)\endalign $$
Now, consider the simple algebraic identity
$$\sum\limits _{i=1}^{n-1} (R_p(a)S_p (a')- S_p (a)R_p(a'))=
 X(a,a')-X(a',a)$$
where
$$X(a,a')=
a^{-1}\tau ^{-1}{\prod (a\tau-b_j)\over
a\tau -a'\tau ^{-1} }-
a^{-1}\tau {\prod (a\tau ^{-1}-b_j)\over
a\tau ^{-1} -a'\tau}  $$
Now let us multiply this identity by
$$\prod\limits _{j=1}^{2n}\tilde{\varphi}(\a ,\b _j)
\tilde{\varphi}(\a' ,\b _j)
L(e^{\a})M(e^{\a '})
e^{{2\pi\over\xi}(\a +\a ')}$$
and integrate the result over $\a$ and $\a '$. The nontrivial thing is
hidden in the integration of the RHS.
Indeed, consider the first part of the corresponding integral:
$$\align
&\int\limits _{-\infty}^{\infty} \int\limits _{-\infty}^{\infty}
\prod\limits _{j=1}^{2n}
\tilde{\varphi}(\a ,\b _j) \tilde{\varphi}(\a',\b _j)
L(e^{\a})M(e^{\a '})X(e^{2\pi \a\over\xi},e^{{2\pi\over\xi}\a '})
e^{{2\pi\over\xi}(\a+\a')}d\a d\a'=\\&=
\int\limits _{-\infty}^{\infty}
\prod\limits _{j=1}^{2n}
\tilde{\varphi} (\a',\b _j)
e^{{2\pi\over\xi}\a '}
M(e^{\a '})       \tag {10}
\\& \left\{\int\limits _{-\infty +2\pi i}^{\infty +2\pi i}
-\int\limits _{-\infty}^{\infty}\right\}
\prod\limits _{j=1}^{2n}
\tilde{\varphi}(\a ,\b _j)
(e^{{2\pi\over\xi}\a }
\tau ^{-1} -e^{{2\pi\over\xi}\b _j})
{e^{-{\pi\over\xi}(\a+\a ')}\over \text{sh}{\pi\over\xi}(\a-\a '-\pi i)}
L(e^{\a})d\a d\a'\endalign $$
The last integral over $\a$ is easy to take since it is sitting on
the simple poles of $1/ \text{sh}{\pi\over\xi}(\a-\a '-\pi i)$. The
pole at $\a =\a '+\pi i$ always exists. For $\xi <\pi$ the additional
poles at $\a =\a '+\pi i\pm i\xi m$ occur but the contributions
from them can be shown to cancel each others. So, the only
contribution comes from the first pole. With the help of (1) the
integral (10) is shown to be expressed in terms of quasiconstants only:
$$\int\limits _{0}^{\infty}\prod\limits _{j=1}^{2n}
{M(A)L(-A)\over A^2+B _j^2}  A^{-1} dA$$
Gathering all pieces we get the deformed version of the Riemann
bilinear identity:
$$\sum\limits _{p=1}^{n-1}
\left(\langle  R_p(a),L(A)\rangle\langle  S_p(a),M(A)\rangle
-\langle  S_p(a),L(A)\rangle\langle  R_p(a),M(A)\rangle \right)
=L\circ M $$
where the intersection number for the polynomials $L$ and $M$ is
defined as
$$ L\circ M=\int\limits _{-\infty}^{\infty}
{L(A)M(-A)-L(-A)M(A) \over
\prod (A^2+B _j^2)}A^{-1}
dA $$
Obviously, the only nonzero intersections are those of even degrees of
$A$ with odd ones. So, the deformed analog of canonic basis can be
constructed in such a way that $A^{2k}$ are the analogues of $a$-cycles
while the $b$-cycles are taken as suitable
linear combinations (with quasiconstant coefficients) of $A^{2k-1}$.
In terms of the canonic basis one can write down the deformed version
of the original form of Riemann bilinear identity.

To finish this paper let us make several comments.
First, as it has been already said, the results obtained here
are needed for the solution of an important physical problem,
namely, the problem of description of all local fields
in SG theory in terms of the form factors bootstrap.
The duality between fields and particles
in SG theory can be described as the duality between deformed
cycles and differentials described here.
So, this
paper is not just a mathematical exercise.
Another interesting point is a possibility of further
deformation. It must be possible to introduce one more
parameter of deformation in such a way that the duality
between polynomials described here is replaced by the duality between
algebraic expressions written in terms of Jacobi $\theta$-functions,
the first step in this direction is done in [5].
\head
Appendix.
\endhead
In this Appendix we explain the regularization of the integral (3)
introduced in [1,2]. Consider the integral
$$\int\limits _{-\infty}^{\infty}\prod\limits _{j=1}^{2n}
\tilde{\varphi} (\a,\b _j)
Q(e^{{2\pi\over\xi}\a })L(e^{\a})
e^{{2\pi\over\xi}\a}d\a $$
which is defined for $\text{deg}(Q(a))\le n-1 $ (the condition
$\text{deg}(L(A))\le 2n-1 $ is always implied).
The definition of the regularized integral for arbitrary polynomial
$Q$ is given by
$$\align
&\int\limits _{-\infty}^{\infty}\prod\limits _{j=1}^{2n}
\tilde{\varphi} (\a,\b _j)
Q(e^{{2\pi\over\xi}\a })L(e^{\a})
e^{{2\pi\over\xi}\a}d\a\equiv
\\&
\int\limits _{-\infty +{\pi i\over 2}-i0}^{\infty +{\pi i\over 2}-i0}
\prod\limits _{j=1}^{2n}
{\tilde{\varphi} (\a,\b _j)\over e^{{2\pi\over\xi}\a }
\tau ^{-3}-e^{{2\pi\over\xi}\b _j}}
Q _1(e^{{2\pi\over\xi}\a })L(e^{\a})
e^{{2\pi\over\xi}\a}d\a +
\\&+
\int\limits _{\Gamma}\prod\limits _{j=1}^{2n}
\tilde{\varphi} (\a,\b _j)
Q _2(\tau ^4 e^{{2\pi\over\xi}\a })L(e^{\a})
e^{{2\pi\over\xi}\a}d\a  \tag {A}
\endalign $$
where the contour $\Gamma $ is drawn around the points
$\a =\b _j -{\pi i\over 2}-i\xi k $ for $k=0,\cdots ,
\left[{\pi\over\xi}\right] $,
polynomials $Q_1$ and $Q_2$ are defined by the
equation:
$$ \align
&Q(a)\prod (a\tau ^{-3}-b_j)=Q_1 (a) +\\&+
 \tau ^{-4}
Q_2(a)\prod (a\tau ^{-1}-b_j) - Q_2(a\tau ^4)\prod (a\tau ^{-3}-b_j)
\endalign $$
It can be shown that the polynomials satisfying this relation and
the additional requirement
$\text{deg}(Q_1)\le 3n-1$
(which is enough for integrals in (A) to converge) can be
found for any $Q$, and that possible ambiguity in definition of
$Q_1,Q_2$ is irrelevant for the value of RHS of (A) [2].

\end